\documentclass[10pt]{amsart}
\usepackage[pdftex]{graphicx}

\usepackage{epsfig}
\newtheorem{thm}{Theorem}
\newtheorem{cor}[thm]{Corollary}
\newtheorem{lem}[thm]{Lemma}
\newtheorem{sat}[thm]{Proposition}
\theoremstyle{definition}

\newtheorem{remark}{Remark}

\usepackage{natbib}
\bibpunct{(}{)}{;}{a}{,}{,}

\usepackage{srcltx}
\usepackage{amssymb,color}

\definecolor{c10}{rgb}{0.1,0.5,0.4}
\definecolor{c15}{rgb}{0.6,0.1,0.4}
\definecolor{c20}{rgb}{0.,0.7,0.}
\definecolor{c30}{rgb}{0.,0.,1.}
\definecolor{c40}{rgb}{1,0.1,0.7}
\definecolor{c50}{rgb}{1,0,0}

\def\cee#1{\textcolor{c20}{#1}}
\def\cee#1{#1}
\newcommand{\nwc}{\newcommand}
\nwc{\COM}[1]{}

  \def\Nset{\mathbb{N}}
  
  \def\Rset{\mathbb{R}}

\newcommand{\neprop}[1]{{Proposition \ref{#1}}}
\newcommand{\netheo}[1]{{Theorem \ref{#1}}}

\newcommand{\E}[1]{\mbox{\rm$ \mathbb{E }$}\left(#1\right)}

\newcommand{\pk}[1]{\mbox{\rm$ \mathbb{P} $}\left(#1\right) }

\def\R{\Rset}
\def\N{\Nset}

\newcommand{\inn}{\in \N}

\newcommand{\QED}{\hfill $\Box$}
\newcommand{\IF}{\infty}

\newcommand{\prooftheo}[1]{ \textbf{Proof of Theorem \ref{#1}}: }
\newcommand{\proofprop}[1]{\textbf{Proof of Proposition \ref{#1}}:}

\newcommand{\BQN}{\begin{eqnarray}}
\newcommand{\EQN}{\end{eqnarray}}
\newcommand{\BQNY}{\begin{eqnarray*}}
\newcommand{\EQNY}{\end{eqnarray*}}
\newcommand{\BS}{\begin{sat}}
\newcommand{\ES}{\end{sat}}
\newcommand{\BL}{\begin{lem}}
\newcommand{\EL}{\end{lem}}
\newcommand{\BT}{\begin{thm}}
\newcommand{\ET}{\end{thm}}
\newcommand{\BK}{\begin{cor}}
\newcommand{\EK}{\end{cor}}
\newcommand{\BD}{\begin{de}}
\newcommand{\ED}{\end{de}}
\newcommand{\BIT}{\begin{itemize}}
\newcommand{\EIT}{\end{itemize}}
\newcommand{\BDI}{\begin{description}}
\newcommand{\EDI}{\end{description}}

\newcommand{\BEX}{\begin{exxa}}
\newcommand{\EEX}{\end{exxa}}

\allowdisplaybreaks[4]

\def\widU{\widetilde{U}}
\def\widhU{\widehat{U}}

\def\LT{\left}
\def\RT{\right}

\begin{document}

\title[Joint density in perturbed classical risk models]{A note on ruin problems in perturbed classical risk models}

\author{Peng Liu}
\address{Peng Liu,  Mathematical Institute, University of Wroc\l aw, pl. Grunwaldzki 2/4, 50-384 Wroc\l aw, Poland, and Department of Actuarial Science, 
University of Lausanne,
UNIL-Dorigny, 1015 Lausanne, Switzerland
}
\email{liupnankaimath@163.com}
\author{Chunsheng  Zhang}
\address{Chunsheng  Zhang,  School of Mathematical Sciences and LPMC, Nankai University, Tianjin 300071, China}
\email{zhangcs@nankai.edu.cn}

\author{Lanpeng Ji}
\address{Lanpeng Ji, Department of Actuarial Science, 
University of Lausanne\\
UNIL-Dorigny, 1015 Lausanne, Switzerland
}
\email{jilanpeng@126.com}

\subjclass[2000]{}

\date{\today}

\begin{abstract}
In this short note, we derive     explicit formulas for the joint densities of the time to ruin and the number of claims until ruin  in perturbed classical risk models, by
constructing several auxiliary random processes.

\end{abstract}

\maketitle

\section{Introduction}

Consider the surplus process $\{U_u(t), t\ge0\}$ of an insurance company modeled by
\begin{eqnarray}\label{model}
U_u(t) = u + ct -\sum_{i=1}^{N(t)}X_i +\sigma B(t),\ \ \ t\geq 0,
\end{eqnarray}
where $u\geq 0$ is the initial reserve, $c>0$ is the rate of premium income, $\{N(t), t\geq 0\}$ is a Poisson claim-counting process with parameter $\lambda>0$, $\{X_i, i\ge1\}$ is a sequence of  independent and identically distributed (iid) positive random variables
with common survival function $\overline{P}(x)$ and density function $p(x)$, representing the
amounts of the successive claims,
 $\{B(t), t\geq 0\}$ is a standard Brownian motion and $\sigma>0$ is a constant representing the diffusion volatility parameter.
 In addition, we assume that $\{N(t), t\ge0\}$, $\{B(t), t\ge0\}$ and $\{X_i, i\ge 1\}$ are mutually independent, and further suppose that the positive net profit condition $c >\lambda \E{X_1}$ holds. In the literature the risk model \eqref{model} is refereed to as the (diffusion) perturbed classical risk model.

 Since the introduction of the  perturbed classical risk model in the seminal contributions \citet{Geber70} and \citet{DufGer1991}, the study of it  has become  popular  in risk theory; see, e.g.,  
 \citet{ZhangWang03}, \citet{Li06} and 
  \citet{Tsai10}. We refer to \citet{AsmAlb2010} for a  nice review  on this model and its generalizations.

Define 
\begin{equation}\label{eq:ruin}
T_u=\inf\{t\ge0: U_u(t)<0\} \ \ (\text{with}\ \inf\{\emptyset\}=\IF)
\end{equation}
to be the time to ruin of the risk model (\ref{model}) with  initial reserve $u\ge0$.  Then $N(T_u)$ is the number of claims until ruin (including the claim which causes ruin).
 In this paper, we are interested in the calculation of the joint density  of $N(T_u)$ and $T_u$  defined by
$$
\omega_u(n,t)=\frac{\partial}{\partial t}\psi_u(n,t),\ \ \ n\in \mathbb{N}_0:=\mathbb{N}\cup\{0\},\ t>0
$$
with
$$\psi_u(n,t)=\pk{N(T_u)=n, T_u\leq t}, \ \ \ n\in\mathbb{N}_0,\ t>0.$$
A nice formula for $\omega_u(n,t)$ in the classical risk model without perturbation (i.e., $\sigma=0$) has been discussed in  \citet{Dick12}. For other related results we refer to \citet{Eg02},  \citet{BorDick2008}, \citet{LanShiWil11}, \citet{Lee11}, \citet{FroPitPol12}, \citet{ZhaoZha13}, \citet{Cheung13}, \citet{LiHJ13}, \citet{LanShi14} and \citet{LiLuJin15}. It is noted that in most of the aforementioned papers results are obtained based on  an application of the known Lagrange's Expansion Theorem. 
For instance, in the framework of \citet{Dick12} (see Eq. (3) therein) the derivation of $\omega_0(n,t)$ relies heavily   on an adequate form of the  inverse Laplace transform of the function $\exp(-\eta_\delta x), \delta, x>0,$ with respect to $\delta$ (which was found by Lagrange's Expansion Theorem). Here $\eta_\delta$ is the unique positive solution of the  Lundberg fundamental
equation
$$
cs -(\lambda+\delta)+ \lambda r \hat{p}(s)=0
$$
with $r\in(0,1]$, $\delta>0$ and $\hat{p}$ the Laplace transform of $p$. We note in passing  that the above  Lundberg  fundamental equation was derived in \citet{LanShiWil11} by introducing the following quantity
\BQN\label{eq:phi}
\phi_{r,\delta}(u)=\E{r^{N(T_u)} e^{-\delta T_u}I_{(T_u<\IF)}}
\EQN
with $I_{(\cdot)}$ the indicator function.

 However, it seems that the approach of  \citet{Dick12} does not  work  well anymore for the perturbed risk model  \eqref{model}. One reason is that now $\omega_0(n,t)=0$ for all $n\inn$ by properties of Brownian motion. In this paper, we shall adopt a different approach.
Similarly to the unperturbed case,  the following generalized Lundberg  fundamental equation 
\BQN\label{eq:rho}
f(s)=\frac{\sigma^2}{2} s^2+cs -(\lambda+\delta)+ \lambda r \hat{p}(s)=0
\EQN
shall play an important role.
 \cee{
In fact, the generalized Lundberg  fundamental equation  was introduced to ensure that, for such $s$ satisfying $f(s)=0$,
$$\{r^{N(t)}e^{ s U_u(t)-\delta t }, {t\geq 0}\}= \{ e^{ s U_u(t)-\delta t +\sum_{i=1}^{N(t)} \ln r }, {t\geq 0}\}$$ 
is a martingale with respect to the natural filtration $\{\mathcal{F}_t^{U_0}, t\geq 0\}$ generated by $\{U_0(t),t\ge0\}$.
Since, for any fixed $s>0$, $\{s U_u(t)-\delta t +\sum_{i=1}^{N(t)} \ln r , t\ge 0\}$ is a spectrally negative L\'{e}vy process, the above martingale property is equivalent to that
$$
\E{e^{ s U_u(t)-\delta t +\sum_{i=1}^{N(t)} \ln r } \big| U_u(0)=u}= e^{su},\ \ \forall t\ge0,
$$
which is valid if and only if $f(s)=0$. 
Note that similar argumentations can also be found in \citet{LanShi14}.
Furthermore, since
$$
f(0)<0, \ \ \lim_{s\rightarrow\IF}f(s)=\IF, \   \ f'(s)=\sigma^2s+c+\lambda r\hat{p}'(s)\geq \sigma^2s+c-\lambda r\E{X_1}\geq \sigma^2s>0,\ \forall s>0,
$$
we conclude that the generalized Lundberg  fundamental equation $f(s)=0$ has a unique positive solution, denoted by $\rho_\delta$. As a result
\BQN \label{eq:MarU}
\{ e^{\widhU_u(t) }, {t\geq 0}\}
\EQN
 with $\widhU_u(t):=\rho_\delta U_u(t)-\delta t +\sum_{i=1}^{N(t)} \ln r, t\ge0$, is a martingale.
}

{\it Outline of this paper:}
 In \neprop{Properho} we derive  an adequate expansion for the key function $\exp(-\rho_\delta x)$ 
 by using a probabilistic approach  instead of resorting to  Lagrange's Expansion Theorem.
With the aid of this expansion and the duality of the risk process $\{U_u(t),t\ge0\}$, we  derive (in \netheo{Thm1}) explicit formulas for the joint density $\omega_u(n,t)$ for any $u>0$. In comparison to the classical analytic approach (as in  \citet{Dick12}) by utilizing Laplace transforms, 
our probabilistic approach  yields much neater formulas. 
The proofs of all the results are relegated to Section 3.



\section{  Results}

Before giving the principle result of this contribution we shall present some preliminary results, among which \neprop{Properho} shall play a crucial role. 
We begin with some further  notation.  For any fixed $x> u$ we define
$$
\tau_u^x=\inf\{t\ge0: U_u(t)=x\}
$$
to be the first hitting time to level $x$ of the risk process $\{U_u(t),t\ge0\}$. Further,
define, for any $t>0,x\in\R$
 \BQN
 &&g_t(0,x) = \frac{1}{2\sqrt{\pi Dt}}e^{-(\lambda t+\frac{(x-ct)^2}{4Dt})}\nonumber\\
&& g_t(n,x) = e^{-\lambda t}\frac{(\lambda t)^n}{n!}\int_{0}^\infty\frac{1}{2\sqrt{\pi Dt}}e^{-\frac{(z+x-ct)^2}{4Dt}}p^{n*}(z)dz,\ \ n\inn,\label{eq:ggt}
 \EQN
 where $D=\sigma^2/2$ and $p^{n*}$, $n\inn$ denotes the $n$-fold convolution of $p$ with itself. It can be shown that $g_t(n,x)$ is the joint density of $N(t)$ and $U_u(t)-u$, \cee{which is independent of the ruin problem considered in the paper.}

 \def \rhod{\rho_\delta}

 \BS \label{Properho}
 Let $\rhod$ be the   unique positive solution of the generalized Lundberg  fundamental equation \eqref{eq:rho}. Then, for any $x> 0$
 \BQN\label{eq:exprho}
 e^{-\rhod x}=\sum_{n=0}^\IF r^n \int_0^\IF e^{-\delta s} \frac{x}{s}g_s(n,x)ds.
 \EQN
  \ES
  \begin{remark}
  {\it The last result  was first derived in \citet{LanShi14} by using an analytic approach. In comparison to the approach therein, our probabilistic approach results in a  much shorter proof.
  }
  \end{remark}
As an application of \neprop{Properho} we obtain  the following  result concerning the joint density of the first hitting time 
and the number of claims until this hitting time. This result is important for the derivation of \neprop{PropHH} given below. 
\BS \label{Prophit}
Let $\tau_u^x$ be the first hitting time to a level $x(>u)$ of the risk process $\{U_u(t),t\ge0\}$. Then
\BQN\label{eq:hitx}
\pk{ N(\tau_u^x)=n,\tau_u^x\in[t,t+dt]}=\frac{x-u}{t}g_t(n, x-u) dt, \ \ \ n\in \N_0,\ t>0.
\EQN

\ES

In order to proceed with the joint density  $\omega_u(n,t)$ for any $u>0$, we  introduce a quantity  $ H(n,t,u,x)$, $n\in  \N_0, t, u, x>0$ as follows (recall $T_u$ is the time to ruin given in \eqref{eq:ruin}) 
  \BQN \label{eq:H1}
   H(n,t,u,x)dx=\pk{N(t)=n,t<T_u,U_u(t)\in[x,x+dx]}.
   \EQN
 Note that, for any fixed $u>0$, 
 $H(n,t,u,x), n\in \mathbb{N}_0,  x>0$ can be interpreted as the joint density of $N(t)$ and $U_u(t)$, with $t$ some fixed time before ruin occurs.
 As it will be seen from our principle result below (\netheo{Thm1}) $H(n,t,u,x)$ is a crucial quantity for  the joint density $\omega_u(n,t)$. We present below an explicit expression for it.
 \BS \label{PropHH}
 Let $ H(n,t,u,x)$, $n\in \mathbb{N}_0, t, u, x>0$ be defined  above. Then
 \BQN \label{eq:HH}
  H(n,t,u,x)=g_t(n,x-u)-\sum_{l=0}^n\int_{0}^t\frac{x}{s} g_s(l,x)g_{t-s}(n-l,-u)ds
 \EQN
 holds for any $n\in \mathbb{N}_0, t, u, x>0$.
 \ES

Now we are ready to present our principle result concerning the joint density  $\omega_u(n,t)$ of $N(T_u)$ and $T_u$ for any $u>0$. It is known that ruin of the perturbed classical risk model \eqref{model} is caused either by a claim or by oscillation; see, e.g., \citet{DufGer1991}.  In the following denote by $\omega_u^s(n,t), n\in\mathbb{N}_0, t>0$  the joint density  when the ruin is caused by a claim, and by $\omega_u^d(n,t), n\in \mathbb{N}_0, t>0$ the joint density when the ruin is caused by oscillation. That is,
$$
\omega_u^s(n,t)=\frac{\partial}{\partial t}\psi_u^s(n,t),\ \ \omega_u^d(n,t)=\frac{\partial}{\partial t}\psi_u^d(n,t), \ n\in\mathbb{N}_0,\ t>0
$$
with
\BQNY
&&\psi_u^s(n,t)=\pk{N(T_u)=n, T_u\leq t, U_u(T_u)<0},\\
&& \psi_u^d(n,t)=\pk{N(T_u)=n, T_u\leq t, U_u(T_u)=0}, \ n\in\mathbb{N}_0,\ t>0.
\EQNY
Clearly, for any $u>0$
\BQN\label{eq:omegaa}
\omega_u(n,t)=\omega_u^s(n,t)+\omega_u^d(n,t),\ \  n \in \mathbb{N}_0, t>0.
\EQN
Our main result below presents  explicit expressions for $\omega_u^s(n,t)$ and $\omega_u^d(n,t)$, $n \in \mathbb{N}_0, t>0,$ and thus in view of the above formula yields an explicit expression for the joint density  $\omega_u(n,t)$, $n \in \mathbb{N}_0, t>0$.  
\BT\label{Thm1}
Let $\omega_u^s(n,t)$ and $\omega_u^d(n,t)$, $n\in \mathbb{N}_0, t>0$  be the joint densities of $N(T_u)$ and $T_u$  defined  above. Then, for any $u>0$
\BQNY
\omega_u^s(0,t)=0,\ \ \ \omega_u^d(0,t)=\frac{u}{2\sqrt{\pi Dt^3}}e^{-(\lambda t+\frac{(u+ct)^2}{4Dt})}, \ \ \ t>0
\EQNY
and
\BQNY
&&\omega_u^s(n,t)=\lambda\int_0^\IF H(n-1,t,u,x)\overline{P}(x)dx, \ \ \ n\in  \N, t>0,\\
&&\omega_u^d(n,t)=\lambda\int_0^t\int_0^\IF\int_0^y H(n-1,t-s,u,y)\frac{y-z}{2\sqrt{\pi Ds^3}}e^{-(\lambda s+\frac{(y-z+cs)^2}{4Ds})}p(z)dzdyds, \ \ \ n\in  \N, t>0,
\EQNY
where $H(n,t,u,x), n\in \mathbb{N}_0, t, u, x>0$  is given in \neprop{PropHH}.
\ET
\COM{
\begin{remark}\label{rem2}
{\it The formulas given in \neprop{PropHH} and \netheo{Thm1} are nice from mathematical point of view, and may be quit useful for simulation purpose. A drawback of them is that closed-form expressions (with reduced number of integrals) can not be easily obtained in general. 
}

\end{remark}
}

\section{Proofs}
In this section we present all the proofs of the results.

\proofprop{Properho} For any fixed $r\in(0,1]$ define an auxiliary process $\{\widU_u(t), t\ge0\}$ from \eqref{model} as follows:
\BQN\label{eq:widUU}
\widU_u(t)=  u + ct -\sum_{i=1}^{N(t)}\widetilde{X}_i +\sigma B(t),\ \ \ t\geq 0,
\EQN
where $\{\widetilde{X}_i,i\ge1\}$ is a sequence of iid generalized positive random variables such that $\widetilde{X}_1$ has density $\widetilde{p}(x)=r p(x), x>0$ and $\pk{\widetilde{X}_1=\IF}=1-r$. It is worth remarking at this point that  the theorems and corollaries in Chapter VII in \citet{Bertoin96} still hold  for  spectrally negative L\'{e}vy processes when the domain of the corresponding L\'{e}vy measure is generalized from $(-\IF,0)$ to $(-\IF,0)\cup\{-\IF\}$. Therefore, by denoting
$$
\widetilde{\tau}_u^x=\inf\{t\ge0: \widetilde{U}_u(t)=x\},\ \ x>u,
$$
we have from Corollary 3 in Chapter VII in \citet{Bertoin96} that the measures $t\pk{\widetilde{\tau}_u^x\in[t,t+dt]}dx$ and $(x-u)\pk{\widetilde{U}_u(t)\in[x,x+dx]}dt$ coincide on $[0,\IF)\times[0,\IF)$. This implies that
\BQN\label{eq:widtau}
\pk{\widetilde{\tau}_u^x\in[t,t+dt]}=\frac{x-u}{t}\widetilde{g}_t(r,x-u)dt,
\EQN
where
\BQNY
\widetilde{g}_t(r,x)&=&\frac{\partial}{\partial x}\pk{\widetilde{U}_u(t)\le x+u}=\frac{1}{2\sqrt{\pi Dt}}e^{-(\lambda t+\frac{(x-ct)^2}{4Dt})}\\
&& \ + \sum_{n=1}^\IF e^{-\lambda t}\frac{(\lambda t)^n}{n!}\int_{0}^\infty\frac{1}{2\sqrt{\pi Dt}}e^{-\frac{(z+x-ct)^2}{4Dt}}(rp)^{n*}(z)dz\\
&&=\sum_{n=0}^\IF r^n g_t(n,x),\ \ x>0.
\EQNY
Moreover, using similar arguments as in Theorem 1 in Chapter VII in \citet{Bertoin96} we obtain that, for any $\delta>0$
\BQN\label{eq:martau}
\E{e^{-\delta \widetilde{\tau}_u^x}}=\E{e^{-\delta \widetilde{\tau}_u^x}I_{(\widetilde{\tau}_u^x<\IF)}}=e^{-\rhod(x-u)},\ \ \ x>u.
\EQN
 Consequently, we conclude from \eqref{eq:widtau}--\eqref{eq:martau} that
\BQNY
 e^{-\rhod (x-u)}=\sum_{n=0}^\IF r^n \int_0^\IF e^{-\delta t} \frac{(x-u)}{t}g_t(n,x-u)dt,\ \ \ x>u,
\EQNY
implying \eqref{eq:exprho}. This completes the proof. \QED

\proofprop{Prophit} 
Recalling \eqref{eq:MarU}, we have that, for any fixed $r\in(0,1], \delta>0$,
 $$\{e^{\widhU_u(t)}, {t\ge 0}\}$$
  is a martingale.
Since $\tau_u^x$ is a stopping time with respect to the filtration $\{\mathcal{F}_t^{U_0}, t\geq 0\}$ and 
$$
e^{\widhU_u(t \wedge \tau_u^x)}\le  e^{\rho_\delta x }<\IF,\ \ \forall t\ge0,
$$
we have from the optional sampling theorem  that 
\BQNY
\E{e^{-\delta \tau_u^x} r^{N(\tau_u^x)} e^{\rhod x}}=\E{e^{\widhU_u(\tau_u^x)}\lvert \widhU_u(0)=\rhod u}=e^{\rhod u}
\EQNY
which means
\BQNY
 \E{e^{-\delta \tau_u^x} r^{N(\tau_u^x)}}= e^{-\rhod (x-u)}.
\EQNY
Consequently, the claim follows by inserting \eqref{eq:exprho} into the last formula. The proof is complete. \QED

\proofprop{PropHH} We introduce the following quantity:
 \BQN \label{eq:H2}
 H^*(n,t,u,x)dx=\pk{N^*(t)=n,t<\tau^{x+u}_{*u},U_u^*(t)\in[x,x+dx]},\ \  n\in \mathbb{N}_0,  t,u,x>0,
\EQN
where, for the fixed $t>0$,
\BQNY
U^*_u(s)=\left\{ \begin{array}{ll}u+U_u(t)-U_u((t-s)-),&0\leq s< t,\\
  U_u(s),&s\geq t,
  \end{array}\right.
\EQNY
$\tau_{*u}^{x+u}$ is defined as
$$
\tau_{*u}^{ x+u}=\inf\{s\ge0: U_u^*(s)=u+x\}, x>0,
$$
and $N^*(t)$ is the number of jumps of the process $\{U_u^*(s),s\ge0\}$ until time $t$. See Figure 1 for the sample paths of $U_u(s)$ and $U_u^*(s)$ when $\sigma=0$.
 \begin{figure}
 \vspace{0mm}
  \includegraphics[width=120mm, height=50mm]{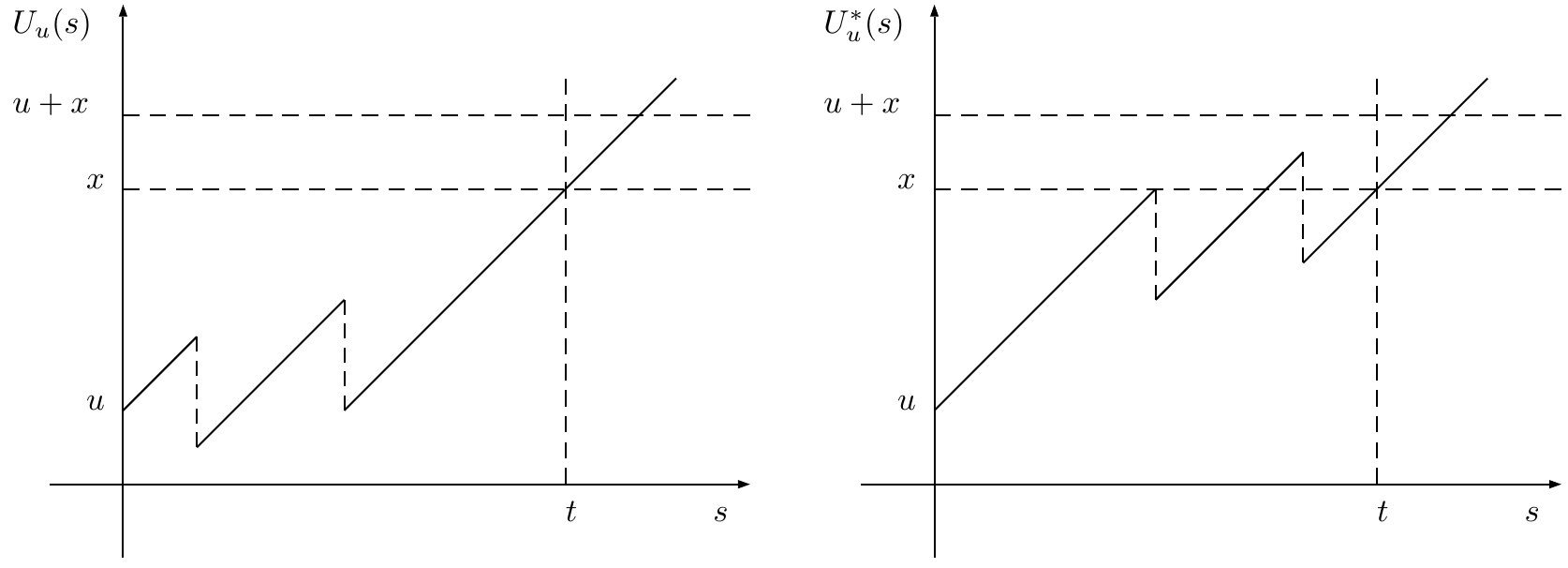}
  \caption{Sample paths of $U_u(s)$ and $U_u^*(s)$ when $\sigma=0$}
\end{figure}
It is noted that the construction of $\{U_u^*(s),s\ge0\}$ from $\{U_u(s),s\ge0\}$ is measure-preserved; see, e.g., \citet{GerShi97} and \citet{GerShi1998}. Therefore, we have
 \BQNY
 \pk{N(t)=n,t<T_u,U_u(t)\in[x,x+dx]}= \pk{N^*(t)=n,t<\tau^{x+u}_{*u},U_u^*(t)\in[x,x+dx]}, 
\EQNY
which means that \BQNY\label{eq:HHH}
 H(n,t,u,x)= H^*(n,t,u,x),\ \  n\in \mathbb{N}_0, t,u, x>0.
 \EQNY
Furthermore, since $\{U^*_u(s), s\geq 0\}$ is also a L\'evy process and $U^*_u(t)=U_u(t)$, by  Theorem 7.10 in \cite{Sato99} 
 we have the process $\{U_u^*(s),  s\ge 0\}$ has the same probability law as the process $\{U_u(s),   s\ge 0\}$. Thus,
\BQNY
H(n,t,u,x)dx&=& H^*(n,t,u,x)dx\\
&=&\pk{N(t)=n,t<\tau^{x+u}_u, U_u(t)\in[x,x+dx]}\\
&=&\pk{N(t)=n, U_u(t)\in[x,x+dx]}- \pk{N(t)=n,t\ge\tau_u^{x+u},U_u(t)\in[x,x+dx]}.
\EQNY
 By the homogeneity and strong Markov property of $\{U_u(s),s\ge0\}$ we conclude that
\BQNY
&&\pk{N(t)=n,t\ge\tau_u^{x+u},U_u(t)\in[x,x+dx]}\\
&&= \sum_{l=0}^n \E{I_{\LT(N(\tau_u^{x+u})=l, \tau_u^{x+u}\le t\RT)} \pk{N(t)-N(\tau_u^{x+u})=n-l, U_{u}(t)\in[x,x+dx]\big\lvert  \mathcal{F}^{U_0}_{\tau_u^{x+u}}}}\\
&&= \sum_{l=0}^n \int_0^t \pk{N(s)=l, \tau_0^{x}\in[s,s+ds]} \pk{N(t-s)=n-l, U_0(t-s)\in[-u,-u+dx]}.
\EQNY
Consequently, the claim  follows from \eqref{eq:ggt} and \eqref{eq:hitx}. This completes the proof. \QED

\prooftheo{Thm1} It follows easily that $\omega_u^s(0,t)=0$. Further, we have
\BQNY
\omega_u^d(0,t)dt&=&\pk{N(t)=0, \inf_{s\in[0,t)}(u+cs+\sigma B(s))>0, \inf_{s\in[t,t+dt]}(u+cs+\sigma B(s))\le 0}\\
&=&\frac{u}{2\sqrt{\pi Dt^3}}e^{-(\lambda t+\frac{(u+ct)^2}{4Dt})} dt,
\EQNY
where  we used the formula for the density of hitting time of a drifted Brownian motion, see, e.g., \cee{Remark 8.3 in Chapter 2 in} \cite{KarShr88}.

Next we consider $\omega_u^s(n,t), n\in\N, t\ge0$. It is noted that $\omega_u^s(n,t)dt$ can be seen as the probability that there are $n-1$ claims until a pre-ruin time $t$ and $U_u(t)\in[x,x+dx]$ and there is a claim (which causes ruin) in $[t, t+dt]$. Thus,
\BQNY
\omega_u^s(n,t)dt&=&\int_0^\IF\int_0^\IF H(n-1,t,u,x)(\lambda dt) p(x+y)dydx\\
&=&\lambda\int_0^\IF H(n-1,t,u,x)  \overline{P}(x )dxdt.
\EQNY
Similarly, $\omega_u^d(n,t)dt$ can be seen as the probability that there are $n-1$ claims  until some pre-last claim occurring time $s$ and  there is a claim  (which does not cause ruin)  in $[s, s+ds]$ and further ruin occurs in  $[t,t+dt]$ by oscillation. Therefore,
\BQNY
\omega_u^d(n,t)dt =  \int_0^\IF  \int_0^y\int_0^t H(n-1,s,u,y)\lambda\ ds\ p(z)\ dz \frac{y-z}{2\sqrt{\pi D(t-s)^3}}e^{-(\lambda (t-s)+\frac{(y-z+c(t-s))^2}{4D(t-s)})}\ dy,
\EQNY
and thus the proof is complete. \QED

\bigskip
{\bf Acknowledgements}:  We are thankful to the associate editor and two referees for several suggestions which significantly improved our manuscript. All the authors kindly acknowledge partial support  by the RARE -318984 (an FP7 Marie Curie IRSES Fellowship)  and Swiss National Science Foundation Project 200021-140633/1. C. Zhang also acknowledges partial support by the National  Science Foundation of China  11371020.

\bibliographystyle{plainnat}
\bibliography{Arch}
\end{document}